\documentclass[10pt]{article}

\usepackage{amsfonts, amssymb, latexsym, amscd}
\usepackage{amsmath, amsthm}
\newtheorem{thm}{Theorem}[section]
\newtheorem{cor}[thm]{Corollary}
\newtheorem{lem}[thm]{Lemma}
\newtheorem{prp}[thm]{Proposition}

\def\1{{\bf 1}}
\def\0{{\bf 0}}
\def\s{{\bf s}}
\def\t{{\bf t}}
\def\uv{{\bf u}}
\def\v{{\bf v}}
\def\w{{\bf w}}
\def\P{{\cal P}}
\def\!{{\hspace{-6pt}}}
\def\ss{\scriptsize}

\title{The graphs with all but two eigenvalues equal to $-2$ or $0$}

\author
{Sebastian M. Cioab{\u{a}}
\\
{\it \small Department of Mathematical Sciences,}
\\
{\it \small University of Delaware, USA. }
\\[3pt]
Willem H. Haemers\thanks{Corresponding author; e-mail haemers@uvt.nl}
\\
{\it\small Department of Econometrics and Operations Research,}
\\
{\it\small Tilburg University, The Netherlands}
\\[3pt]
Jason R. Vermette
\\
{\it \small Natural Sciences Division,}
\\
{\it \small Missouri Baptist University, Saint Louis, USA}
\\[10pt]
\noindent
Dedicated to Andries E. Brouwer on the occasion of his 65th birthday.
}

\date{}

\begin{document}

\maketitle

\begin{abstract}
\noindent
We determine all graphs for which the adjacency matrix has at most
two eigenvalues (multiplicities included) not equal to $-2$, or $0$,
and determine which of these graphs are determined by their
adjacency spectrum.
\end{abstract}

\section{Introduction}
In an earlier paper~\cite{chvw} we determined the class $\cal G$ of
graphs with at most two adjacency eigenvalues
(multiplicities included) different from $\pm 1$.
The classification was motivated by the question whether the
friendship graph is determined by its spectrum.
Here we deal with the class $\cal H$ of graphs for which the
adjacency matrix $A$ has all but at most two eigenvalues equal
to $-2$ or $0$.
Equivalently, $A+I$ has at most two eigenvalues different from
$\pm 1$.
The class $\cal H$ is in some sense complementary to $\cal G$.
Indeed, many graphs in $\cal G$ (including the friendship graphs)
are the complements of graphs in $\cal H$.
In particular it easily follows that a regular graph is in $\cal G$
if and only if its complement is in $\cal H$.
But for non-regular graphs there is no such relation.

Note that a graph from $\cal H$ remains in $\cal H$ if we add or
delete isolated vertices.
Let $\cal H'$ be the set of graphs in $\cal H$ with no isolated
vertices.
Then it clearly suffices to determine all graphs in $\cal H'$.

It turns out that $\cal H'$ contains several infinite families and some sporadic graphs,
and that relatively few graphs in $\cal H'$ have a cospectral mate.
Thus we find many graphs with four distinct eigenvalues which are determined by their adjacency spectrum;
see~\cite{dkx} for some discussion on this phenomenon.

The major part of our proof deals with graphs in $\cal H'$ with two positive eigenvalues.
These graphs have least eigenvalues $-2$ and have been classified, see for example~\cite{crs}.
However, using this classification is not straightforward, and we decided to start from scratch
using mainly linear algebra.
The most important tool is {\em eigenvalue interlacing},
which states that if $\lambda_1(A)\geq\cdots\geq\lambda_n(A)$
are the eigenvalues of a symmetric matrix $A$ of order $n$,
and $\lambda_1(B)\geq\cdots\geq\lambda_m(B)$
are the eigenvalues of a principal submatrix $B$ of $A$ or order $m$,
then
\[
\lambda_i(A)\geq\lambda_i(B)\geq\lambda_{n-m+i}(A)
\ \mbox{ for }\ i=1,\ldots,m.
\]
Another important tool is the following well known result on
equitable partitions (we refer to~\cite{bh} for these tools and
many other results on spectral graph theory).
Consider a partition ${\P}=\{V_1,\ldots,V_m\}$ of the set
$V=\{1,\ldots,n\}$.
The characteristic matrix $\chi_\P$ of $\P$ is the $n\times m$
matrix whose columns are the character vectors of $V_1,\ldots,V_m$.
Consider a symmetric matrix $A$ of order $n$, with rows and
columns partitioned according to $\P$.
The partition of $A$ is {\em equitable} if each submatrix
$A_{i,j}$ formed by the rows of $V_i$ and the columns of $V_j$
has constant row sums $q_{i,j}$.
The $m\times m$ matrix $Q=(q_{i,j})$ is called the
{\em quotient matrix} of $A$ with respect to $\P$.
\begin{lem}\label{ep}
The matrix $A$ has the following two kinds of eigenvectors
and eigenvalues:
\begin{itemize}
\item[(i)]
The eigenvectors in the column space of $\chi_\P$;
the corresponding eigenvalues coincide with the eigenvalues of $Q$.
\item[(ii)]
The eigenvectors orthogonal to the columns of $\chi_\P$;
the corresponding eigenvalues of $A$ remain unchanged if
some scalar multiple of the all-one block $J$ is added to block
$A_{i,j}$ for each $i,j\in\{1,\ldots,m\}$.
\end{itemize}
\end{lem}
The all-ones and all-zeros vector of dimension $m$ is denoted by
$\1_m$ (or $\1$) and $\0_m$ (or $\0$), respectively.
We denote the $m\times n$ all-ones matrix by $J_{m,n}$ (or just $J$),
the all-zeros matrix by $O$, and the identity matrix of order $n$ by
$I_n$, or $I$.

\section
{Graphs in $\cal{H'}$ with just one positive eigenvalue}\label{1}

It is known (see~\cite{s}) that a graph with just one positive
adjacency eigenvalue is a complete multipartite graph,
possibly extended with a number of isolated vertices.
Thus we have to find the complete multipartite graphs in $\cal H'$.
\begin{thm}\label{onepos}
If a graph $G\in{\cal H'}$ has one positive eigenvalue,
then $G$ is one of the following.
\begin{itemize}
\item $G_0(\ell,m)=K_{\ell,m}$ with spectrum
$\{0^{\ell+m-2},~\pm\sqrt{\ell m}\}$ ($\ell\geq m\geq 1$),
\item $G_1=K_{1,1,3}$ with spectrum $\{-2,~-1,~0^2,~3\}$,
\item $G_2(k,m)=K_{2,\ldots,2,m}$ with
spectrum $\{-2^{k-2},~0^{k+m-2},~k-2\pm\sqrt{k^2+2(k-1)(m-2)}\}$
($m\geq 1$), where $k\geq 3$ is the number of classes.
\end{itemize}
\end{thm}

\noindent
{\bf Proof.}
Let ${\cal H}_1'$ be the set of complete multipartite graph in
$\cal H'$, and let $G$ be a complete $k$-partite graph in
${\cal H}_1'$.
Then $k\geq 2$, and $k=2$ gives $G_0(\ell,m)$.
Assume $k\geq 3$.
The graph $K_{2,3,3}$ has two eigenvalues less than $-2$,
so by interlacing,
$G$ does not contain $K_{2,3,3}$ as an induced subgraph.
So, if no class has size $1$, then $G=G_2(k,m)$ with
$m\geq 2$.

Suppose there are $k_1\geq 1$ classes of size~$1$ and assume
$k\geq k_1+3$.
Then at least one class has size $2$, and at most one class has
size greater than $2$
(because $K_{2,3,3}$ is a forbidden induced subgraph).
If one class has size $m\geq 3$, then the adjacency matrix $A$
of $G$ admits an equitable partition with quotient matrix
\[
Q=\left[
\begin{array}{ccc}
k_1-1 & m & 2(k-k_1-1)\\
k_1   & 0 & 2(k-k_1-1)\\
k_1   & m & 2(k-k_1-2)
\end{array}
\right].
\]
By Lemma~\ref{ep}, the eigenvalues of $Q$ are also eigenvalues of $A$,
but $Q$ has three eigenvalues different from $-2$ and $0$,
so $G\not\in{\cal H}_1'$.
If $k_1\geq 3$, then $A+I$ has at least $3$ equal rows,
and hence $A$ has an eigenvalue $-1$ of multiplicity at least $2$,
and therefore $G\not\in{\cal H}'_1$.
If $k_1=2$ and $G$ contains $K_{2,3}$ (which has smallest eigenvalue
$-\sqrt{6}<-2$), then $G$ has an eigenvalue equal to $-1$,
and an eigenvalue less than $-2$, so $G\not\in{\cal H}_1'$.

So only $K_{1,2,\ldots,2}$, $K_{1,\ell,m}$, $K_{1,1,m}$,
and $K_{1,1,2,\ldots,2}$ are not covered yet by the above cases,
and one easily checks that only $K_{1,2,\ldots,2}=G_2(k,1)$ and $K_{1,1,3}=G_1$ survive.
The eigenvalues readily follow by use of Lemma~\ref{ep},
or from~\cite{eh}.
\hfill $\Box$
\\

The graph $K_{2,\ldots,2}$ with $k\geq 1$ classes is also known as
the {\em cocktail party graph}, usually denoted by $CP(k)$.
We will write $C(k)$ for the adjacency matrix of $CP(k)$.
The spectrum of $CP(k)$ equals $\{-2^{k-1},~0^k,~2k-2\}$.
Note that $CP(k)$ with $k\geq 2$ and $K_{1,4}$ are the only graphs
in $\cal H'$ with just one eigenvalue different from $-2$ or $0$.
Clearly there is no graph in $\cal H'$ with no positive eigenvalue,
so we have the following result:
(The disjoint union of two graphs $G$ and $G'$ is denoted by $G+G'$,
and $mG$ denotes $m$ disjoint copies of $G$.)

\begin{thm}\label{one}
Suppose $G$ is a graph with $n$ vertices and at most one
eigenvalue different from $-2$ and $0$,
then $G=nK_1$, $G=K_{1,4}+(n-5)K_1$, or $G=CP(k)+(n-2k)K_1$ with
$k\geq 2$.
\end{thm}

\begin{cor}\label{disconnected}
If $G$ is a disconnected graph in $\cal H'$, then $G$ is one of
the following:
\begin{itemize}
\item
$G_3=2K_{1,4}$ with spectrum $\{-2^2,~0^{8},~2^2\}$,
\item
$G_4(k)=K_{1,4}+CP(k)$ with spectrum $\{-2^k,~0^{k+3},~2,~2k-2\}$
($k\geq 2$),
\item
$G_5(k,\ell)=CP(k)+CP(\ell)$ with spectrum
$\{-2^{k+\ell-2},~0^{k+\ell},~2\ell-2,~2k-2\}$ ($k\geq\ell\geq 2$).
\end{itemize}
\end{cor}
\noindent
{\bf Proof.}
Suppose $G\in{\cal H'}$ has $c$ components, then each component has
at least one edge, and therefore at least one positive
eigenvalue.
Thus $G$ has $c$ positive eigenvalues.
So, if $G$ is disconnected then $c=2$, and each component is
$K_{1,4}$ or $CP(k)$ with $k\geq 2$.
\hfill$\Box$\\

\section
{Graphs in $\cal H'$ with two positive eigenvalues}\label{2}

In this section we determine the set ${\cal H}_2'$ of connected
graphs in ${\cal H}'$ with two positive eigenvalues.
We start with some tools.

\begin{prp}\label{E}
Let $A$ be the adjacency matrix of a graph $G\in{\cal H}_2'$,
and define $E=A(A+2I)$.
Then $E$ is positive semi-definite and rank$(E)=2$.
\end{prp}

We call a pair of vectors $\v$ and $\w\in\{0,1\}^n$ {\em almost equal}
whenever $\v$ and $\w$ have different weights, and differ in at
most two positions.

\begin{cor}\label{ae}
If $A$ is the adjacency matrix of a graph $G\in{\cal H}_2'$,
then no two columns of $A+I$ are equal or almost equal.
\end{cor}
\noindent{\bf Proof.}
Let $a\geq b$ be the weights of two (almost) equal columns of $A+I$.
Then $a-b\leq 2$ and $E=A(A+2I)=(A+I)^2-I$ has the following
principal submatrix
\[
\left[
\begin{array}{cc}
a-1& b\\b & b-1
\end{array}
\right],
\]
which is not positive semi-definite.
\hfill$\Box$
\\

\begin{prp}
The following graphs ($(d)$-$(m)$ represented by their adjacency matrices)
are forbidden induced subgraphs for any graph $G\in{\cal H}_2'$.
\[
(a)\ K_{1,5}
,
\ \ (b)\ K_{2,3}
,
\ \ (c)\ C_5
,
\ \ (d)
{\ss
\left[
\begin{array}{ccccc}
0&\! 0&\! 1&\! 1&\! 1\\
0&\! 0&\! 1&\! 1&\! 1\\
1&\! 1&\! 0&\! 0&\! 0\\
1&\! 1&\! 0&\! 0&\! 1\\
1&\! 1&\! 0&\! 1&\! 0
\end{array}
\right]
},
\ (e)
{\ss
\left[
\begin{array}{ccccc}
0&\! 0&\! 0&\! 1&\! 1\\
0&\! 0&\! 0&\! 1&\! 1\\
0&\! 0&\! 0&\! 0&\! 1\\
1&\! 1&\! 0&\! 0&\! 0\\
1&\! 1&\! 1&\! 0&\! 0
\end{array}
\right]
},
\]
\[
\ (f)
{\ss
\left[
\begin{array}{cccccc}
0&\! 0&\! 0&\! 1&\! 1&\! 1\\
0&\! 0&\! 0&\! 1&\! 1&\! 1\\
0&\! 0&\! 0&\! 1&\! 1&\! 1\\
1&\! 1&\! 1&\! 0&\! 1&\! 1\\
1&\! 1&\! 1&\! 1&\! 0&\! 1\\
1&\! 1&\! 1&\! 1&\! 1&\! 0
\end{array}
\right]
},
\ (g)
{\ss
\left[
\begin{array}{cccccc}
0&\! 0&\! 0&\! 1&\! 1&\! 1\\
0&\! 0&\! 0&\! 0&\! 1&\! 1\\
0&\! 0&\! 0&\! 0&\! 1&\! 1\\
1&\! 0&\! 0&\! 0&\! 1&\! 1\\
1&\! 1&\! 1&\! 1&\! 0&\! 1\\
1&\! 1&\! 1&\! 1&\! 1&\! 0
\end{array}
\right]
},
\ (h)
{\ss\left[
\begin{array}{cccccc}
0&\! 0&\! 0&\! 1&\! 1&\! 1\\
0&\! 0&\! 0&\! 0&\! 1&\! 1\\
0&\! 0&\! 0&\! 0&\! 1&\! 1\\
1&\! 0&\! 0&\! 0&\! 0&\! 0\\
1&\! 1&\! 1&\! 0&\! 0&\! 1\\
1&\! 1&\! 1&\! 0&\! 1&\! 0
\end{array}
\right]
},
\,(i)\,\,
{\ss
\left[
\begin{array}{cccccc}
0&\! 0&\! 0&\! 0&\! 1&\! 1\\
0&\! 0&\! 0&\! 0&\! 1&\! 1\\
0&\! 0&\! 0&\! 0&\! 1&\! 1\\
0&\! 0&\! 0&\! 0&\! 1&\! 1\\
1&\! 1&\! 1&\! 1&\! 0&\! 1\\
1&\! 1&\! 1&\! 1&\! 1&\! 0
\end{array}
\right]
},
\]
\[
\ (j)
{\ss\left[
\begin{array}{cccccc}
0&\! 0&\! 0&\! 0&\! 1&\! 1\\
0&\! 0&\! 0&\! 0&\! 1&\! 1\\
0&\! 0&\! 0&\! 0&\! 1&\! 1\\
0&\! 0&\! 0&\! 0&\! 0&\! 1\\
1&\! 1&\! 1&\! 0&\! 0&\! 1\\
1&\! 1&\! 1&\! 1&\! 1&\! 0
\end{array}
\right]
},
\ (k)
{\ss
\left[
\begin{array}{cccccc}
0&\! 0&\! 0&\! 0&\! 1&\! 1\\
0&\! 0&\! 0&\! 0&\! 0&\! 1\\
0&\! 0&\! 0&\! 0&\! 0&\! 1\\
0&\! 0&\! 0&\! 0&\! 0&\! 1\\
1&\! 0&\! 0&\! 0&\! 0&\! 1\\
1&\! 1&\! 1&\! 1&\! 1&\! 0
\end{array}
\right]
},
\ (\ell)
{\ss
\left[
\begin{array}{cccccc}
0&\! 0&\! 0&\! 0&\! 1&\! 1\\
0&\! 0&\! 0&\! 0&\! 0&\! 1\\
0&\! 0&\! 0&\! 0&\! 0&\! 1\\
0&\! 0&\! 0&\! 0&\! 0&\! 1\\
1&\! 0&\! 0&\! 0&\! 0&\! 0\\
1&\! 1&\! 1&\! 1&\! 0&\! 0
\end{array}
\right]
},
\ (m)
{\ss
\left[
\begin{array}{cccccc}
0&\! 0&\! 0&\! 0&\! 1&\! 1\\
0&\! 0&\! 0&\! 0&\! 1&\! 1\\
0&\! 0&\! 0&\! 0&\! 0&\! 1\\
0&\! 0&\! 0&\! 0&\! 0&\! 1\\
1&\! 1&\! 0&\! 0&\! 0&\! 0\\
1&\! 1&\! 1&\! 1&\! 0&\! 0
\end{array}
\right]
}.
\]

\end{prp}
\noindent{\bf Proof.}
Interlacing.
Each of the adjacency matrices has smallest eigenvalue less than
$-2$, or more than two positive eigenvalues.
\hfill$\Box$
\\

Let $\alpha=\alpha(G)$ be the maximum size of a coclique (that is,
an independent set of vertices) in $G$.
We will distinguish a number of cases depending on $\alpha$.

\subsection{Graphs in ${\cal H}_2'$ with $\alpha=2$}

\begin{thm}\label{bipartitecomplement}
If $G\in{\cal H}_2'$ and $\alpha(G)=2$,
then ${G}$ is one of the following graphs (represented by their adjacency matrices):
\begin{itemize}
\item
$G_6(m):~
\left[
\begin{array}{cc}
J-I_m & I_m \\
I_m & J-I_m
\end{array}
\right]
$
($m\geq 3$), with spectrum
$\{-2^{m-1},\, 0^{m-1},\, m-2,\, m\}$,
\item
$G_7:~
\left[
\begin{array}{ccc}
J-I_7 & J-I_7 & \0 \\
J-I_7 & J-I_7 & \1 \\
\0^\top & \1^\top & 0
\end{array}
\right]
$
with spectrum $\{-2^7,\, 0^6,\, 7\pm 2\sqrt{7}\}$.
\end{itemize}
\end{thm}
\noindent
{\bf Proof.}
Let $\overline{G}$ be the complement of $G$.
Since $\alpha(G)=2$, $\overline{G}$ has no triangles.
It is well-known that the complement of an odd cycle $C_n$
with $n\geq 5$ has at least three positive eigenvalues.
Therefore, by interlacing, $\overline{G}$ has no induced
odd cycle, so $\overline{G}$ is bipartite and the
adjacency matrix $A$ of $G$ has the following structure:
\[
A=\left[
\begin{array}{cc}
J-I_m & N\\
N^\top & J-I_{m'}
\end{array}
\right].
\]
First we claim that $|m-m'|\leq 1$.
Indeed, suppose $m\leq m'-2$, then $J_{m'}-I$ has eigenvalue $-1$ of multiplicity $m'-1$,
and therefore, by interlacing, $A$ has an eigenvalue $-1$ of multiplicity at least $m'-1-m>0$,
contradiction.
So without loss of generality $m=m'\geq 2$ or $m=m'-1\geq 2$.

Consider four columns $\s$, $\t$, $\uv$, $\v$ of $A+I$, with
$\s$ and $\t$ in the first part and $\uv$ and $\v$ in the second part.
Let $m+a$, $m+b$ ($a\leq b$), $m'+c$, and $m'+d$ ($c\leq d$)
be the weights of $\s$, $\t$, $\uv$ and $\v$ respectively,
and define $\lambda=\s^\top \t -m$ and $\mu=\uv^\top \v -m'$.
Then these four columns give the following submatrix of $E=A(A+2I)$:
\[
E'=\left[
\begin{array}{cccc}
m+a-1 & m+\lambda & a+c & a+d \\
m+\lambda & m+b-1 & b+c & b+d \\
a+c & b+c & m'+c-1 & m'+\mu \\
a+d & b+d & m'+\mu & m'+d-1
\end{array}
\right].
\]
After Gaussian elimination (subtract the first row from second row
and the last row from the third row,
and apply a similar action to the columns) we obtain
\[
E''=\left[
\begin{array}{cccc}
m+a-1 & \lambda-a+1 & c-d & a+d \\
\lambda-a+1 & a+b-2\lambda-2 & 0 & b-a \\
c-d & 0 & c+d-2\mu-2 & \mu-d+1 \\
a+d & b-a & \mu-d+1 & m'+d-1
\end{array}
\right].
\]
By Lemma~\ref{E}, rank$(E'')\leq 2$, which implies that the upper
$3\times 3$ submatrix of $E''$ has determinant $0$, which leads to
\[
(m+a-1)(a+b-2\lambda-2)(c+d-2\mu-2)=(c+d-2\mu-2)
(\lambda-a+1)^2+(a+b-2\lambda-2)(c-d)^2.
\]
By Lemma~\ref{ae}, no two rows of $A+I$ are equal or almost equal
which implies that
$a+b-2\lambda-2\geq 0$ with equality if and only if $a=b=\lambda+1$,
and similarly
$c+d-2\mu-2\geq 0$ with equality if and only if $c=d=\mu+1$.
If $a+b-2\lambda-2$ and $c+d-2\mu-2$ are both nonzero, then we get
\[
m+a-1=(a-\lambda-1)^2 /(a+b-2\lambda-2) + (d-c)^2/(c+d-2\mu-2).
\]
By use of $a\leq b$, $c\leq d$ and $d\geq\mu+1$ we find
$m+a-1\leq(a-\lambda-1)/2 + (d-c)/2 < m$, contradiction.
Therefore $a=b=\lambda+1$ or $c=d=\mu+1$.
When $m=m'$ this gives $N=I_m$ or $N=J-I_m$.
If $m=2$ or $N=J-I_m$, then $G=CP(m)$, and so only $N=I_m$ with
$m\geq 3$ survives.
When $m=m'-1$ we find four possibilities:
$N=[~I_m\ \ \0~],\ [~J-I_m\ \ \1~],\ [~I_m\ \ \1~]$,
or $[~J-I_m\ \ \0~]$.
The first two options do not occur, because $A+I$ has almost
equal columns (but note that the second case gives
$K_{1,2,\ldots,2}$).
For the other cases $A$ has equitable partitions with
quotient matrices
\[
\left[
\begin{array}{ccc}
m-1 & 1   & 1 \\
1   & m-1 & 1 \\
m   & m   & 0
\end{array}
\right]
\mbox{ and }
\left[
\begin{array}{ccc}
m-1 & m-1 & 0 \\
m-1 & m-1 & 1 \\
0   & m & 0
\end{array}
\right],
\]
respectively.
These quotient matrices, and therefore the adjacency matrices, have
three eigenvalues different from $-2$ and $0$, unless $m=2$ in the
first case (which gives $K_{1,2,2}$), and $m=7$ in the second case.
\hfill$\Box$
\\

\subsection{Graphs in ${\cal H}_2'$ with $\alpha\geq 3$}

\begin{thm}\label{remaining}
If $G\in{\cal H}_2'$ and $\alpha(G)\geq 3$, then $G$ is one of
the following graphs (represented by their adjacency matrices):
\begin{itemize}
\item
$G_8(k,\ell):~
\left[
\begin{array}{ccc}
C(k) & O & \1 \\
O & \! C(\ell)\! & \1 \\
\1^\top & \1^\top & 0
\end{array}
\right]$
with spectrum
$\{-2^{k+\ell-1},\ 0^{k+\ell},\ k+\ell-1\pm\sqrt{(k-\ell)^2 + 1}\}$
\\[3pt]
($k\geq\ell\geq 1$, $k\geq 2$),
\item
$G_9(k):~
\left[
\begin{array}{cccc}
O & O & J_{2,3} & J_{2,2k} \\
O & J-I & J-I & O \\
J_{3,2} & J-I & J-I & J_{3,2k} \\
J_{2k,2} & O & J_{2k,3} & C(k)
\end{array}
\right]$
with spectrum $\{ -2^{k+3},\, 0^{k+3},\, k+3\pm\sqrt{k^2+3}\}$
\\[3pt]
($k\geq 0$, if $k=0$, the last block row and column vanish),
\item
$G_{10}:~$
{\scriptsize
$\left[
\begin{array}{ccccccccc}
0\!&0\!&0\!&0\!&1\!&1\!&1\!&1\!&0\\
0\!&0\!&0\!&0\!&0\!&0\!&0\!&0\!&1\\
0\!&0\!&0\!&0\!&0\!&0\!&0\!&0\!&1\\
0\!&0\!&0\!&0\!&1\!&1\!&1\!&1\!&1\\
1\!&0\!&0\!&1\!&0\!&0\!&1\!&1\!&0\\
1\!&0\!&0\!&1\!&0\!&0\!&1\!&1\!&1\\
1\!&0\!&0\!&1\!&1\!&1\!&0\!&0\!&0\\
1\!&0\!&0\!&1\!&1\!&1\!&0\!&0\!&1\\
0\!&1\!&1\!&1\!&0\!&1\!&0\!&1\!&0
\end{array}
\right]
$
}
with spectrum $\{-2^3,\, 0^4,\, 3\pm\sqrt{2}\}$,
\item
$G_{11}:~$
{\scriptsize
$\left[
\begin{array}{ccccccccc}
0\!&0\!&0\!&0\!&0\!&1\!&1\!&0\!&0\\
0\!&0\!&0\!&0\!&0\!&1\!&1\!&0\!&0\\
0\!&0\!&0\!&0\!&0\!&0\!&0\!&1\!&1\\
0\!&0\!&0\!&0\!&0\!&0\!&0\!&1\!&1\\
0\!&0\!&0\!&0\!&0\!&1\!&1\!&1\!&1\\
1\!&1\!&0\!&0\!&1\!&0\!&1\!&0\!&1\\
1\!&1\!&0\!&0\!&1\!&1\!&0\!&1\!&0\\
0\!&0\!&1\!&1\!&1\!&0\!&1\!&0\!&1\\
0\!&0\!&1\!&1\!&1\!&1\!&0\!&1\!&0
\end{array}
\right]
$
} with spectrum $\{-2^3,\ 0^4,\ 2,\ 4\}$,
\item
$G_{12}:~$
{\scriptsize
$\left[
\begin{array}{cccccccc}
0\!&0\!&0\!&0\!&0\!&0\!&1\!&0\\
0\!&0\!&0\!&0\!&0\!&0\!&1\!&0\\
0\!&0\!&0\!&0\!&0\!&0\!&0\!&1\\
0\!&0\!&0\!&0\!&0\!&0\!&0\!&1\\
0\!&0\!&0\!&0\!&0\!&0\!&1\!&1\\
0\!&0\!&0\!&0\!&0\!&0\!&1\!&1\\
1\!&1\!&0\!&0\!&1\!&1\!&0\!&1\\
0\!&0\!&1\!&1\!&1\!&1\!&1\!&0
\end{array}
\right]
$
} with spectrum $\{-2^2,\ 0^4,\ 1,\ 3\}$.
\end{itemize}
\end{thm}

\noindent{\bf Proof.}
Let $C$ be a coclique of $G$ of size $\alpha\geq 3$.
We assume that $C$ has the largest number of outgoing edges among all
cocliques of size $\alpha$.
The coclique $C$, and the remaining vertices of $G$ give the
following partition of the adjacency matrix $A$ of $G$.
\[
A=\left[
\begin{array}{cc}
O & N \\ N^\top & B
\end{array}
\right]
\]
First we will prove that the matrix $N$ defined above has one of
the following structures:
\[
(i)\ \left[
\begin{array}{cc}
J_{p,a} & O \\
O       & J_{q,b}
\end{array}
\right],
\ (ii)\ \left[
\begin{array}{ccc}
J_{p,a} & O      & J_{p,c}\\
O       & J_{q,b} & J_{q,c}
\end{array}
\right],
\ (iii)\ \ \left[
\begin{array}{cc}
J_{p,a} & J_{p,c} \\
O       & J_{q,c}
\end{array}
\right],
\ (iv)\ \left[
\begin{array}{cc}
J_{p,a} & O \\
O       & J_{q,b}\\
J_{r,a} & J_{r,b}
\end{array}
\right].
\]
Suppose rank$(N)\geq 3$.
Then $N$ has a $3\times 3$ submatrix of rank $3$, and $A$ has
a nonsingular principal submatrix $A'$ of order $6$ containing
$O$ of order $3$.
Interlacing gives $\lambda_3(A')\geq \lambda_3(O) = 0$.
But $A'$ is nonsingular, so $\lambda_3(A')>0$, and interlacing
gives $\lambda_3(G)>0$, contradiction.
So rank$(N)\leq 2$.
Suppose $N=J$, then $B=J-I$, since otherwise $G$ contains the
forbidden subgraph $K_{2,3}$.
Thus we have $G=K_{\alpha,1,\ldots,1}$
which has been treated in Section~\ref{1}.
We easily have that $N$ has no all-zero column
(because $C$ is maximal),
and no all-zero row (because $G$ is connected).
We conclude that rank$(N)=2$ and $N$ has one of the structures
$(i),(ii),(iii),(iv)$.
We will continue the proof by considering these four cases step
by step.

\subsubsection{Case $(i)$}\label{i}
We have
\[
A=\left[\begin{array}{cccc}
O&O&J_{p,a} & O \\
O&O& O      & J_{q,b}\\
J_{a,p}& O & B_1 & M \\
O & J_{b,q} & M^\top & B_2
\end{array}
\right]
\]
Assume $p\leq q$.
If $p=1$, then $B_1=J-I$ (because $C$ is maximal), and $M=O$
(otherwise there would be another coclique of size $\alpha$
with more outgoing edges).
So $G$ is disconnected, contradiction.
If $p\geq 3$, then the forbidden subgraphs $(a),(b),(f)$, and $(i)$
lead to just six possibilities for $(p,q,a,b)$, being:
\[
(3,3,2,2),\ (3,3,2,1),\ (3,3,1,1),\ (3,4,2,1),\ (3,4,1,1),
\ (4,4,1,1).
\]
One easily checks that none of these give a graph in ${\cal H}_2'$.

So $p=2$.
We assume $a\geq b$ (if $p=q=2$ we can do so, and if $q>p$,
then the same forbidden subgraphs as above imply $b\leq2$,
but if $b=2$ and $a=1$ then the last two columns of $A+I$ are equal
or almost equal, which contradicts Corollary~\ref{ae}).
Forbidden subgraph $(d)$ shows that each row of $B_1$ has weight
$a-1$, or $a-2$.
If $B_1$ has an off-diagonal $0$, then the two corresponding rows
have weight $a-2$, and $G$ has another coclique of size $\alpha$.
The two corresponding rows of $M$ have weight $0$, because otherwise
the other coclique would have more outgoing edges.
So if all rows of $B_1$ have weight $a-2$, then $M=O$,
and $G$ is disconnected.
Therefore at least one row (the first one, say) of $B_1$ has
weight $a-1$.
Now rows 1, 3 and $\alpha+1$ of $A$ give the following submatrix of
$E=A(A+2I)$:
\[
E'=\left[
\begin{array}{ccc}
a & 0 & a+1\\
0 & b & x \\
a+1 & x & 1+a+x
\end{array}
\right],
\]
where $x$ is the weight of the first row of $M$.
By Proposition~\ref{E}, rank$(E)=2$, so $\det(E')=0$,
which leads to $b/a=bx-x^2-b$.
Therefore $b/a$ is an integer, and since $a\geq b$,
we have $a=b$ and $2x=a\pm\sqrt{(a-2)^2-8}$.
This equation only has an integer solution if $a=b=5$ and
$x=2$, or $3$.
Forbidden subgraphs $(b)$ and $(f)$ show that $q<3$, so $p=q=2$.
Now also $B_2$ has a row (the first one say) with weight $b-1=4$,
and the corresponding row of $A$ together with the three rows
considered above, give a $4\times 4$ submatrix of $E$ of rank~$3$,
contradiction with Proposition~\ref{E}.

\subsubsection{Cases $(ii)$ and $(iii)$}\label{iii}
We consider case $(iii)$, but include $(ii)$ by
allowing $a=0$, or $b=0$ (not both).
The forbidden subgraph $K_{1,5}$ gives $p+q=\alpha\leq 4$.
First consider $\alpha=4$.
Then forbidden subgraphs $(b)$ and $(i)$
give $c=1$, and $(b)$, $(k)$, $(\ell)$ and $(j)$ imply that $p=q=2$.
We assume $a\geq b$.
Moreover, the forbidden graph $(m)$ shows that the last vertex of
$A$ is adjacent to all other vertices.
So $A$ has the following structure
\[
A=\left[
\begin{array}{ccccc}
O&O&J_{2,a} & O &\1 \\
O&O& O      & J_{2,b}&\1 \\
J_{a,2}& O & B_1 & M&\1 \\
O & J_{b,2} & M^\top & B_2 &\1 \\
\1^\top & \1^\top &\1^\top &\1^\top &0
\end{array}
\right].
\]
Like in Section~\ref{i}, $p=2$ implies that a row of $B_1$ has
weight $a-1$ or $a-2$.
If some row of $B_1$ (the first one, say) has weight $a-1$, then
row 1, 2 and 5 of $A$ give the following submatrix of $E=A(A+2I)$:
\[
E'=\left[
\begin{array}{ccc}
a+1 & 1 & a+1\\
1 & b+1 & x+1 \\
a+1 & x+1 & 1+a+x
\end{array}
\right],
\]
where $x$ is the weight of the first row of $M$.
By Proposition~\ref{E} $\det(E')=0$, which implies that $x\neq b$
and $x=1+1/(a+1)+1/(b-x)$, which has no solution (because $a\geq b$).
Therefore each row of $B_1$ has weight $a-2$, and $B_1=C(a/2)$.
But then $M=O$, since otherwise $G$ would have another coclique of
size~$4$ with more outgoing edges.
Also $B_2=C(b/2)$ ($CP(0)$ is the graph with the empty vertex set),
because if a row of $B_1$ would have weight $b-1$,
then the corresponding row of $A+I$ and row~3 of $A+I$ are almost
equal, which contradicts Corollary~\ref{ae}.
Now reordering rows and columns of $A$ shows that $G=G_8(k,\ell)$).
The eigenvalues readily follow by use of Lemma~\ref{ep}.
\\

Next we consider $\alpha=3$.
Then we take $p=1$, $q=2$.
The forbidden graphs $(b)$ and $(f)$ show that $c=1$ or $2$.
First suppose $c=1$.
Then
\[
A=\left[
\begin{array}{ccccc}
O&O&\1^\top & O &1 \\
O&O& O      & J_{2,b}&\1 \\
\1& O & B_1 & M&\v \\
O & J_{b,2} & M^\top & B_2 &\w \\
1 & \1^\top &\v^\top &\w^\top &0
\end{array}
\right].
\]
We have $B_1=J-I$, because $C$ is maximal.
Let $x$ be the weight of an arbitrary row of $M$.
If $x=0$, then the corresponding row of $A+I$ is equal or almost equal
to the first row of $A+I$.
So $x\geq 1$.
Because no coclique of size $3$ has more outgoing edges than $C$,
it follows that $x=1$ and $\v=\0$.
Forbidden subgraph $(e)$ gives $\w=\1$.
But now the last row of $A+I$ is almost equal to the second row.
So there is no solution when $\alpha=3$ and $c=1$.
(Note that the conclusion also holds if $a=0$ or $b=0$.)

Suppose $c=2$.
We assume $p=1$, $q=2$.
Then $A$ has the following structure.
\[
A=\left[
\begin{array}{ccccc}
O&O& J_{1,a} & O &\1^\top \\
O&O& O      & J_{2,b}&J_{2,2} \\
J_{a,1} & O & B_1 & M& M_1 \\
O & J_{b,2} & M^\top & B_2 &M_2 \\
\1 & J_{2,2} & M_1^\top &M_2^\top & B_3
\end{array}
\right].
\]
Like before, we have that $B_1=J-I$.
Forbidden subgraphs $(b)$ and $(e)$ give that
$B_3=J-I$, and $M_2=J$.
Note that if $a=0$, then the last two rows of $A+I$
are equal, and if $b=0$ the first and the fourth row are
(almost) equal.
So $a\geq 1$, and $b\geq 1$.
Moreover, subgraphs $(g)$ and $(h)$ imply that
each row of $M_1$ equals $[0~1]$, or $[1~0]$.
Also, by the same arguments as before each row of $M$
has weight $1$, and each row of $B_2$ has weight
$b-1$ or $b-2$.
Consider the submatrix of $E$ made from the second, the fourth
and the last two rows of $A$,
and then replace the last two rows and columns by their sum.
Then we get
\[
E'=\left[
\begin{array}{ccc}
b+2 & 2 & 2b+6 \\
2 & a+2 & a+6 \\
2b+6 & a+6  & a+4b+18
\end{array}
\right].
\]
We have rank$(E')\leq\mbox{rank}(E)=2$,
so $0=\det(E')=2(a-2)(ab+2a+2b)$, hence $a=2$.
Thus $M_1=[\1_2 \ O]$ or $[I_2 \ O]$.
However, the second options does not occur, because it
would give a column in $A+I$ which is almost equal to the last
column of $A+I$.
Similarly, if a column of $M$, different from the first one,
has weight $b-1$, then the corresponding column in $A+I$ is
almost equal to the last one.
So these columns of $M$ all have weight~$b-2$.
However, the first row and column of $M$ has weight $b-1$,
since if some off-diagonal entry (the second, say) of
the first row of $M$ is zero, then rows and columns
2, 3, 4, 6, 7 give the forbidden subgraph $(e)$.
Now there is just one family of solutions for $A$,
which, after reordering rows and columns shows that $G=G_9(k)$.
The spectrum follows in a straightforward way by use of
Lemma~\ref{ep}.
This finishes cases~$(ii)$ and $(iii)$.

\subsubsection{Case ($iv$)}

We have the following adjacency matrix.
\[
A=\left[
\begin{array}{ccccc}
O&O&O& J_{p,a} & O \\
O&O&O& O       & J_{q,b} \\
O&O&O& J_{r,a} & J_{r,b} \\
J_{a,p} & O & J_{a,r} & B_1 & M \\
O & J_{b,q} & J_{b,r} & M^\top & B_2
\end{array}
\right].
\]
We assume $p\leq q$.
There is a somewhat special behaviour when $p=r=b=1$, so
we treat this case first.
Forbidden subgraph $(a)$, $(b)$ and $(d)$ show that $q\leq 3$ and that each row of $B_1$ has weight at least $a-2$.
If some row of $B_1$ has weight $a-1$, then $A+I$ has two almost equal rows, hence $B_1$ is the adjacency matrix of $CP(a/2)$.
Two nonadjacent vertices in $CP(a/2)$ together with $q$ vertices of $C$ make another coclique in $G$ of size $\alpha$ which,
by assumption, has not more outgoing edges than $C$.
Therefore the two corresponding entries of $M$ cannot both be equal to $1$ and therefore the weight $x$ of $M$ satisfies $x\leq a/2$.
The first, the second and the last row of $A$ give the following submatrix of $E=A(A+2I)$:
\[
E'=\left[
\begin{array}{ccc}
a & 0 & x \\
0 & 1 & 2 \\
x & 2 & x+q+1
\end{array}
\right].
\]
Proposition~\ref{E} gives $\det(E')=0$, which implies $2x=a\pm\sqrt{(a+2(q-3))^2-4(q-3)^2}$.
If $q=1$ we get $a=9$, which is not possible.
If $q=2$ we get $a=4$, $x=2$, which gives $G_{10}$.
If $q=3$ there is no solution, because $G$ contains a forbidden subgraph $(k)$ or $(\ell)$.
\\

Next we treat the case $p=q=r=1$ and $b\geq 2$.
Then, like before, forbidden subgraphs $(b)$ and $(d)$ imply
that the a row of $B_1$ has weight $a-1$ or $a-2$.
Suppose all rows have weight $a-2$, so $B_1$ represents $CP(a/2)$.
Now the first, the second and $k$-th row of $A$
($4\leq k\leq a+3$) give the following submatrix of $E=A(A+2I)$:
\[
E'=\left[
\begin{array}{ccc}
a & a & 0 \\
a & a+x & x \\
0 & x & b
\end{array}
\right],
\]
where $x$ is the weight of the corresponding row in $M$.
This matrix is singular, which implies that $x=0$, or $x=b$.
Since $b\geq 2$, no row or column of $B_2$ has weight $b$, since
otherwise $A+I$ would have two equal or almost equal columns.
So $B_2$ is the adjacency matrix of $CP(b/2)$,
and the columns of $M$ have weight $0$ or $a$.
So $M=O$ or $M=J$, but $M=J$ is impossible,
because then one vertex of $C$ and two nonadjacent vertices of
$CP(a/2)$ make up a coclique of size~3
with more outgoing edges than than $C$ has.
Thus we find an equitable partition of $A$ with quotient matrix
\[
Q=\left[
\begin{array}{ccccc}
0&0&0& a & 0 \\
0&0&0& 0 & b \\
0&0&0& a & b \\
1&0&1& a-2 & 0 \\
0&1&1& 0 & b-2
\end{array}
\right].
\]
For all positive integers $a$ and $b$, $Q$ and hence $A$ has at
least three eigenvalues different from $-2$ and $0$.

We conclude that at least one row of $B_1$ has
weight $a-1$ and find the following submatrix $E'$ of $E=A(A+2I)$:
\[
E'=\left[
\begin{array}{ccc}
a & 0 & a+1\\
0 & b & x \\
a+1 & x & a+x+1
\end{array}
\right],
\]
where $x$ is the weight of the corresponding row of $M$.
Now $\det(E')=0$ gives $x^2-xb+b+b/a = 0$,
We assume $a\geq b$.
Then $b=a$ and $2x=b\pm\sqrt{(b-2)^2-8}$.
This only has an integer solution if $b=5$.
Then $x=2$ or $x=3$.
It is straightforward to check that there is
no matrix $A$ that satisfies all requirements.
[A quick way to see this is by extending $E'$ above to a $4\times 4$
matrix $E''$ by considering also a row of $B_2$ with weight~$4$.
If $x'$ is the weight of the corresponding row of $M^\top$,
then rank$(E'')=2$ implies that $x=2$, $x'=3$ (or vice versa).
This is clearly impossible if all rows of $B_1$ and $B_2$ have
weight~$4$.
However, if a number of rows of $B_1$ and $B_2$ have weight $3$,
then this number is even and the corresponding rows of $M$ and
$M^\top$ have weight $0$ or $5$ (see above).
This gives that the sum of row weights of $M$ is equal to
$2$ or $1$ (mod~$5$), whilst the column weights add up to
$3$ or $4$ (mod~$5$),which is of course impossible.]
This finishes the case $p=q=r=1$.

The next case is $p=r=1$, $q\geq 2$.
Forbidden subgraphs $a$ and $b$ give $q=2$ and $b\leq 2$,
or $q=3$, $b=1$.
The case $b=1$ is solved above, so we only need to deal with
$b=q=2$.
Then forbidden subgraph $(b)$ gives that $B_2=J-I$.
We consider the submatrix $E'$ of $E=A(A+2I)$ formed by the first,
the second and the last row of $A$.
Then
\[
E'=\left[
\begin{array}{ccc}
a & 0 & x\\
0 & 2 & 3 \\
x & 3 & x+4
\end{array}
\right],
\]
where $x$ is the weight of the second column of $M$.
We have $0=\det(E')=a(2x-1)-2x^2$, which implies $x=1$, $a=2$.
Now it is straightforward to check that there is no solution.

Finally we consider the case that $3\leq p+r\leq q+r$.
Forbidden subgraphs $(a)$ and $(b)$ lead to $p+r\leq 4$,
$a\leq 2$, $b\leq 2$,
$a\leq 1$ if $p+r=4$, and $b\leq 1$ if $q+r=4$.
This gives only a short list of possibilities for $(p,q,r,a,b)$,
being:
\[
\begin{array}{cccccc}
(2,2,1,2,2),&
(2,2,1,2,1),&
(2,2,1,1,1),&
(2,3,1,2,1),&
(2,3,1,1,1),&
(3,3,1,1,1),
\\
(1,1,2,2,2),&
(1,1,2,2,1),&
(1,1,2,1,1),&
(1,2,2,2,1),&
(1,2,2,1,1),&
(2,2,2,1,1).
\end{array}
\]
Only the first and the last $ 5$-tuple correspond to a solution:
$G_{11}$ and $G_{12}$, respectively.
\hfill$\Box$\\

\section{Conclusions}
The  graphs in ${\cal H}'$ are given in
Theorems~\ref{onepos}, \ref{bipartitecomplement}, \ref{remaining}
and Corollary~\ref{disconnected}.
The graphs $G_8(2,2)$ and $G_{11}$ are the only two nonisomorphic
cospectral graphs in ${\cal H}'$.
All other pairs of cospectral graphs in ${\cal H}$ consist of pairs of
graphs in ${\cal H}'$ with the same nonzero part of the spectrum
(see Theorem~\ref{conclusion} below),
where one or both graphs are extended with isolated vertices,
such that the numbers of vertices in both graphs are equal.
The following theorem follows straightforwardly from the list of graphs in ${\cal H}'$.
\begin{thm}\label{conclusion}
Two nonisomorphic graphs $G,~G'\in{\cal H}'$ have equal nonzero parts of the spectrum if and only if one of the following holds:
\begin{itemize}
\item
$G=G_0(\ell,m)\ (= K_{\ell,m})$, $G'=G_0(\ell',m')\ (= K_{\ell',m'})$, where $\ell m=\ell'm'$,
\item
$G,~G'\in\{G_4(k),~G_5(k,2)\}$ ($k\geq 2$),
\item
$G,~G'\in\{G_5(k+1,k),~G_6(2k),~G_8(k,k)\}$ ($k\geq 2$),
\item
$G,~G'\in\{G_6(3),~G_{12}\}$,
\item
$G,~G'\in\{G_4(4),~G_5(4,2),~G_9(1)\}$,
\item
$G,~G'\in\{G_4(3),~G_5(3,2),~G_6(4),~G_8(2,2),~G_{11}\}$,
\item
$G,~G'\in\{G_3,~G_4(2),~G_5(2,2)\}$.
\end{itemize}
\end{thm}
\begin{cor}
A graph $G\in{\cal H}'$ is not determined by the spectrum of the adjacency matrix if and only if $G$ is one of the following:
\begin{itemize}
\item
$G_0(\ell,m)\ (=K_{\ell,m})$, where $\ell m$ has a divisor
strictly between $\ell$ and $m$,
\item
$G_4(k),~G_5(k+1,k),~G_8(k,k)$ with $k\geq 2$,
\item
$G_3,~G_5(4,2),~G_{11},~G_{12}$.
\end{itemize}
\end{cor}
For every graph in $\cal H$ we can decide whether it is determined
by the spectrum of the adjacency matrix by use of
Theorem~\ref{conclusion}.
For example $CP(k)+ K_1$, which is the complement of the friendship
graph, is determined by its spectrum because every graph with
the same nonzero part of the spectrum equals $CP(k)+mK_1$ for some $m\geq 0$.
This result was already proved in \cite{ajo}.
Also the spectral characterisation of complete multipartite graphs
is known, see~\cite{eh}.

\section*{Acknowledgments} The research of S.M. Cioab\u{a} was supported by the NSF Grant DMS 160078.

\end{document}